\documentclass[]{article}

\usepackage[margin=1.5in]{geometry} %toglierlo

\usepackage{amsmath}
\usepackage{amssymb}
\usepackage{amsthm}

\newtheorem{theorem}{Theorem}

\DeclareMathOperator{\cl}{cl} % per la chiusura di un insieme

\begin{document}

\title{A Functional Analytic Approach For A Singularly Perturbed Dirichlet Problem For The Laplace Operator In A Periodically Perforated Domain.}

\author{Paolo Musolino}

\date{}
\maketitle

\noindent
{\bf Abstract:} We consider a sufficiently regular bounded open connected subset $\Omega$ of $\mathbb{R}^n$ such that $0 \in \Omega$ and such that $\mathbb{R}^n \setminus \cl\Omega$ is connected. Then we choose a point $w \in ]0,1[^n$. If $\epsilon$ is a small positive real number, then we define the periodically perforated domain $T(\epsilon) \equiv \mathbb{R}^n\setminus \cup_{z \in \mathbb{Z}^n}\cl(w+\epsilon \Omega +z)$. For each small positive $\epsilon$, we introduce a particular Dirichlet problem for the Laplace operator in the set $T(\epsilon)$. More precisely, we consider a Dirichlet condition on the boundary of the set $w+\epsilon \Omega$, and we denote the unique periodic solution of this problem by $u[\epsilon]$. Then we show that (suitable restrictions of) $u[\epsilon]$ can be continued real analytically in the parameter $\epsilon$ around $\epsilon=0$. 

\vspace{11pt}



%%%%%%%%%%%%%%%%%%%%%%%%%%%%%%%%%%%%%%%%%%%%
%% MAINMATTER
%%%%%%%%%%%%%%%%%%%%%%%%%%%%%%%%%%%%%%%%%%%%

\section{Introduction}

This paper is devoted to  the analysis of a singularly perturbed Dirichlet problem for the Laplace operator in a periodically perforated domain.

Throughout this paper $n \in \mathbb{N} \setminus \{0,1\}$. Let $\alpha \in ]0,1[$. We consider an open bounded connected subset $\Omega$ of $\mathbb{R}^n$ of class $C^{1,\alpha}$ such that $0 \in \Omega$ and such that the complement in $\mathbb{R}^n$ of the closure $\cl \Omega$ is also connected. Then let $w$ be a point of the fundamental cell $A \equiv ]0,1[^n$. For each $\epsilon \in \mathbb{R} \setminus \{0\}$ we consider the open bounded connected subset $\Omega_\epsilon$ defined by
\[
\Omega_\epsilon \equiv w+\epsilon \Omega.
\]
Then we take $\epsilon_0>0$ such that
\[
\cl \Omega_\epsilon \subseteq A \qquad \forall \epsilon \in ]0,\epsilon_0[.
\]
If $\epsilon \in  ]0,\epsilon_0[$, then the set $\Omega_\epsilon$ represents the perforation of the fundamental cell $A$. If $z \in \mathbb{Z}^n \setminus \{0\}$, we make a hole in the set $A+z$, by considering the translation $\Omega_\epsilon+z$ of the perforation $\Omega_\epsilon$ of the fundamental cell $A$.
Next, if $\epsilon \in  ]0,\epsilon_0[$, we consider the periodically perforated set
\[
T(\epsilon)\equiv \mathbb{R}^n \setminus \cup_{z \in \mathbb{Z}^n}\cl (z+\Omega_\epsilon),
\]
obtained by removing from $\mathbb{R}^n$ the set $\cup_{z \in \mathbb{Z}^n}\cl (z+\Omega_\epsilon)$.

In order to introduce a Dirichlet boundary value problem for the Laplace operator in the perforated set $T(\epsilon)$, we take a function
\[
g \in C^{1,\alpha}(\partial \Omega).
\] 
We note that $g$ is defined on the fixed set $\partial \Omega$, not depending on $\epsilon$. Then, for each $\epsilon \in ]0,\epsilon_0[$, we shall consider the following perturbed Dirichlet problem:
 \begin{equation}\label{bvp:Direps}
 \left \lbrace 
 \begin{array}{ll}
 \Delta u (x)= 0 & \textrm{$\forall x \in T(\epsilon)$}, \\
u(x+e_i) =u(x) &  \textrm{$\forall x \in \mathrm{cl} T(\epsilon), \quad \forall i \in \{1,\dots,n\}$}, \\
u(x)=g\bigl(\frac{1}{\epsilon}(x-w)\bigr) & \textrm{$\forall x \in \partial \Omega_\epsilon$},
 \end{array}
 \right.
 \end{equation}
where $\{e_1,\dots,e_n\}$ is the canonical basis of $\mathbb{R}^n$. If $\epsilon \in ]0,\epsilon_0[$, then problem \eqref{bvp:Direps} admits a unique solution in $C^{1,\alpha}(\cl T(\epsilon))\cap C^2(T(\epsilon))$, and we denote it by $u[\epsilon]$.

Then we pose the following questions:
\begin{enumerate}
\item[(i)] Let $\bar{x}$ be fixed in $\mathbb{R}^n \setminus (w+\mathbb{Z}^n)$. What can be said on the map $\epsilon \mapsto u[\epsilon](\bar{x})$ around $\epsilon=0$ for $\epsilon$ small and positive?
\item[(ii)] What can be said on the map $\epsilon \mapsto \mathcal{E}(\epsilon)\equiv \int_{A\setminus \cl \Omega_\epsilon}|\nabla u[\epsilon](x)|^2 dx$ around $\epsilon=0$ for $\epsilon$ small and positive?
\end{enumerate}

Questions of this type have long been investigated for problems on a bounded domain with a small hole with the methods of asymptotic analysis, which aims at giving complete asymptotic expansions of the solutions in terms of the parameter $\epsilon$. It is perhaps difficult to provide a complete list of the contributions. Here, we mention the work of Kozlov, Maz'ya and Movchan \cite{KoMaMo99}, Maz'ya, Nazarov and Plamenewskij \cite{MaNaPl00}, Ozawa \cite{Oz83}, Vogelius and Volkov \cite{VoVo00}, Ward and Keller \cite{WaKe93}. We also mention the vast literature of homogenization theory (cf. \textit{e.g.}, Dal Maso and Murat \cite{DaMu04}).

Here we wish to characterize the behaviour of $u[\epsilon]$ at $\epsilon=0$ by a different approach. Thus for example, if we consider a certain functional, say $f(\epsilon)$, relative to the solution such as, for example, one of those considered in questions (i)-(ii) above, one could resort to Asymptotic Analysis and may succeed (depending on the functional $f$ under consideration) to write out an expansion of the type
\begin{equation}\label{eq:fas}
f(\epsilon)=\sum_{j=0}^ra_j\epsilon^j +o(\epsilon^r) \qquad \text{as $\epsilon\to 0^+$},
\end{equation}
for suitable coefficients $a_j$. Instead, in the same circumstance we would try to prove that $f(\cdot)$ can be continued real analytically around $\epsilon=0$. More generally, we would try to represent $f(\epsilon)$ for $\epsilon>0$ in terms of real analytic maps and in terms of possibly singular at $\epsilon=0$, but known functions of $\epsilon$ (such as $\epsilon^{-1}$, $\log \epsilon$, etc.). We observe that our approach does have certain advantages. Indeed, if, for example, we know that $f(\epsilon)$ equals for $\epsilon>0$ a real analytic function of $\epsilon$ defined in a neighbourhood of $\epsilon=0$ in the real line, then we know that an asymptotic expansion as \eqref{eq:fas} for all $r$ would necessarily generate a convergent series
\[
\sum_{j=0}^\infty a_j \epsilon^j,
\]
and that the sum of such a series would be $f(\epsilon)$ for $\epsilon>0$.

Such a project has been carried out by Lanza de Cristoforis in many papers for the case of a bounded domain with a small hole. Here we mention in particular Lanza \cite{La08} (see also references therein) where a Dirichlet boundary value problem for the Laplace operator is considered. We also mention Dalla Riva and Lanza \cite{DaLa10}.

As far as problems in periodically perforated domains are concerned, we mention, for instance, the work of Ammari, Kang and Touibi \cite{AmKaTo05}, where a linear transmission problem is considered in order to compute an asymptotic expansion of the  effective electrical conductivity of a periodic dilute composite. Furthermore, we note that periodically perforated domains are extensively studied in the frame of homogenization theory. Among the vast literature, here we mention, \textit{e.g.}, Cioranescu and Murat \cite{CiMu82, CiMu83}, Ansini and Braides \cite{AnBr02}. We also observe that boundary value problems in domains with periodic inclusions can be analysed with the method of functional equations (at least for the two dimensional case). Here we mention, \textit{e.g.}, the work of Mityushev and Adler \cite{MiAd02}, where a doubly periodic Dirichlet problem for the Poisson equation is studied in order to compute the longitudinal permeability of spatially periodic rectangular arrays of circular cylinders. We also mention Rogosin, Dubatovskaya, and Pesetskaya \cite{RoDuPe09}, where Eisenstein functions are used to construct the solution of a mixed boundary value problem for a doubly periodic multiply connected domain, in order to study effective properties of a doubly periodic 2D composite material.

\section{Strategy and main results}

We now briefly outline our strategy. First of all, we observe that problem \eqref{bvp:Direps}, which we consider only for positive $\epsilon$, is singular (and not defined) for $\epsilon=0$. Then, by exploiting a potential theoretic approach, we can convert problem \eqref{bvp:Direps} into an equivalent integral equation. Indeed, if $\epsilon \in ]0,\epsilon_0[$, then the solution $u[\epsilon]$ can be represented as the sum of a periodic double layer potential and a costant. The density of the periodic double layer potential and the costant must solve a particular integral equation defined on the $\epsilon$-dependent domain $\partial \Omega_\epsilon$. We observe that the periodic layer potentials are constructed by replacing the fundamental solution of the Laplace equation $S_n$ with a periodic analogue $S_n^a$  in the definition of the classical layer potentials for the Laplace operator. Such a periodic analogue is a periodic distribution $S_n^a \in \mathcal{D}'(\mathbb{R}^n)$ such that
\[
\Delta S_n^a= \sum_{z \in \mathbb{Z}^n}\delta_z-1,
\]
in the sense of distributions. For a construction of $S_n^a$, we refer, for instance, to Hasimoto \cite{Ha59}, Shcherbina \cite{Sh86}, Poulton, Botten, McPhedran and Movchan \cite{PoBoMcMo99}, Ammari, Kang and Touibi \cite{AmKaTo05}. Then we observe that by changing the variables appropriately, we can obtain an equivalent integral equation which can be analysed by means of the Implicit Function Theorem around the degenerate case in which $\epsilon=0$, and we represent the unknowns of the integral equation in terms of real analytic functions of $\epsilon$ defined in a neighbourhood of $0$. By exploiting these results we can prove our main Theorem.

\begin{theorem}\label{thm:ueps}
Let $V$ be a bounded open subset of $\mathbb{R}^n$, such that 
\[
\cl V \subseteq \mathbb{R}^n \setminus (w+\mathbb{Z}^n).
\]
Then there exist $\epsilon_1 \in ]0,\epsilon_0[$, a real analytic map $U_1$ of $]-\epsilon_1,\epsilon_1[$ to $C^0(\cl V)$ and a real analytic map $U_2$ of $]-\epsilon_1,\epsilon_1[$ to $\mathbb{R}$, such that 
\[
\cl V \subseteq T(\epsilon) \quad \forall \epsilon \in ]0,\epsilon_1[,
\]
and
\[
u[\epsilon](x)=\epsilon^{n-1}U_1[\epsilon](x)+U_2[\epsilon] \qquad \forall x \in \cl V, \quad \forall \epsilon \in ]0,\epsilon_1[.
\]
\end{theorem}

We observe that Theorem \ref{thm:ueps} implies that $u[\epsilon]_{|\cl V}$ converges uniformly to a constant, namely $U_2[0]$, as $\epsilon\to 0^+$. On the other hand, we note that one may expect that $u[\epsilon]$ converges (in some sense) to a solution of some problem on the unperturbed domain $\mathbb{R}^n$. Actually, this is what happens, since constant functions are easily seen to be the solutions of the following problem:
\[
 \left \lbrace 
 \begin{array}{ll}
 \Delta u (x)= 0 & \textrm{$\forall x \in \mathbb{R}^n$}, \\
u(x+e_i) =u(x) &  \textrm{$\forall x \in \mathbb{R}^n, \quad \forall i \in \{1,\dots,n\}$}. \\
  \end{array}
 \right.
\]

As far as the energy integral of the solution is concerned, we have the following.

\begin{theorem}\label{thm:eneps}
There exist $\epsilon_2 \in ]0,\epsilon_0[$ and a real analytic map $G$ of $]-\epsilon_2,\epsilon_2[$ to $\mathbb{R}$, such that 
\[
\int_{A\setminus \cl \Omega_\epsilon}|\nabla u[\epsilon](x)|^2dx=\epsilon^{n-2}G[\epsilon] \qquad \forall \epsilon \in ]0,\epsilon_2[.
\]
\end{theorem}

\section{Conclusion}

We note that, by an approach similar to the one used for problem \eqref{bvp:Direps}, we can study different boundary value problems in the periodically perforated domain $T(\epsilon)$. Moreover, we can also consider the case where the fundamental cell $A$ is of the form $\prod_{i=1}^n]0,\bar{a}_i[$, with $\{\bar{a}_1,\dots,\bar{a}_n\} \subseteq ]0,+\infty[$, and the case where also the length of the sides of the fundamental cell goes to $0$.

\section*{Acknowledgments}
 This paper is an announcement based on the work performed by the author in his ``Laurea Magistrale'' Thesis \cite{Mu08} under the guidance of  M. Lanza de Cristoforis.

\bibliographystyle{aipproc}

\begin{thebibliography}{99}

\bibitem{KoMaMo99} V.~A. Kozlov, V.~G. Maz'ya and A.~B. Movchan, \emph{Asymptotic Analysis of Fields in Multistructures,} New York: The Clarendon Press, Oxford University Press, 1999.

\bibitem{MaNaPl00} V.~G. Maz'ya, S.~A. Nazarov and B.~A. Plamenewskij. \emph{Asymptotic Theory of Elliptic Boundary Value Problems in Singularly Perturbed Domains,} Basel: Birkh\"auser Verlag, 2000.

\bibitem{Oz83} S.~Ozawa, \emph{J. Fac. Sci. Univ. Tokyo Sect. IA Math.} \textbf{30}, 53-62 (1983).

\bibitem{VoVo00} M.~S. Vogelius and D.~Volkov, \emph{M2AN Math. Model. Numer. Anal.} \textbf{34}, 723-748 (2000).

\bibitem{WaKe93} M.~J. Ward and J.~B. Keller, \emph{SIAM J. Appl. Math.} \textbf{53}, 770-798 (1993).

\bibitem{DaMu04} G.~Dal Maso and F.~Murat, \emph{Ann. Inst. H. Poincaré, Anal. Non Linéaire} \textbf{21}, 445-486 (2004).

\bibitem{La08} M.~Lanza de Cristoforis, \emph{Analysis (Munich)} \textbf{28}, 63-93 (2008).

\bibitem{DaLa10} M.~Dalla Riva and M.~Lanza de Crisotoris, \emph{Analysis (Munich)} \textbf{30}, 67-92 (2010).

\bibitem{AmKaTo05}H.~Ammari, H.~Kang and K.~Touibi,  \emph{Asymptot. Anal.} \textbf{41}, 119-140 (2005).

\bibitem{CiMu82} D.~Cioranescu and F.~Murat, ``Un terme étrange venu d'ailleurs'', in  \emph{Nonlinear Partial Differential Equations and Their Applications. Collège de France Seminar}, Vol. II, edited by H.~Brezis and J.~L. Lions, Res. Notes in Math. 60, Pitman, London, 1982, pp. 98-138, 389-390. 

\bibitem{CiMu83} D.~Cioranescu and F.~Murat, ``Un terme étrange venu d'ailleurs, II'', in  \emph{Nonlinear Partial Differential Equations and Their Applications. Collège de France Seminar}, Vol. III, edited by H.~Brezis and J.~L. Lions, Res. Notes in Math. 70, Pitman, London, 1983, pp. 154-178, 425-426.

\bibitem{AnBr02} N.~Ansini and A.~Braides, \emph{J. Math. Pures Appl.}  \textbf{81}, 439-451 (2002).

\bibitem{MiAd02} V.~Mityushev and P.~M. Adler, \emph{Z. Angew. Math. Mech.} \textbf{82}, 335-345 (2002).

\bibitem{RoDuPe09} S.~Rogosin, M.~Dubatovskaya, E.~Pesetskaya, \emph{\v{S}iauliai Math. Sem.} \textbf{4}, 167-187 (2009).

\bibitem{Ha59} H.~Hasimoto, \emph{J. Fluid Mech.} \textbf{5}, 317-328 (1959).

\bibitem{Sh86} V.~A. Shcherbina, \emph{Teor. Funktsi\u{\i} Funktsional. Anal. i Prilozhen} \textbf{45},  132-139 (1986). 

\bibitem{PoBoMcMo99} C.~G. Poulton, L.~C. Botten, R.~C. McPhedran and A.~B. Movchan, \emph{Proc. R. Soc. Lond. A} \textbf{455}, 1107-1123  (1999).

\bibitem{Mu08} P.~Musolino, ``Due problemi di perturbazione singolare su domini con perforazioni multiple. Un approccio funzionale analitico'', Laurea Magistrale Thesis, Universit\`a di Padova, 2008.

\end{thebibliography}

\end{document}